\documentclass{elsarticle}
\usepackage{graphicx}
\usepackage{amsmath,amsthm,amsfonts}

\def\cubic{\mathrm{cubic}}

\newcommand\R{\mathbb{R}}
\newcommand\C{\mathbb{C}}
\newcommand\Z{\mathbb{Z}}

\newcommand\T{\mathbb{T}}
\newcommand\calH{\mathcal{H}}

\newcommand\calB{\mathcal{B}}
\newcommand\calE{\mathcal{E}}

\newcommand\calL{\mathcal{L}}

\newcommand\self{\mathrm{\mathrm{self}}}

\newcommand\dist{{\mathrm{dist}}}
\newcommand\rel{{\mathrm{rel}}}

\newcommand\eps{{\varepsilon}}

\numberwithin{equation}{section}
\begin{document}
\begin{frontmatter}

\title{Dislocation microstructures and  strain-gradient plasticity 
with one active slip plane}


\author[iam]{Sergio Conti\corref{mycorrespondingauthor}}
\author[roma]{Adriana Garroni}
\author[iam]{Stefan M\"uller}
\cortext[mycorrespondingauthor]{Corresponding author.\\
Tel. 0049 228 7362211, Fax 0049 228 73 62251, Email sergio.conti@uni-bonn.de}
\address[iam]{Institut f\"ur Angewandte Mathematik,
Universit\"at Bonn, 53115 Bonn, Germany}
\address[roma]{Dipartimento di Matematica, Sapienza, Universit\`a di Roma,
00185 Roma, Italy}
%\address{Institut f\"ur Angewandte Mathematik,
%Universit\"at Bonn 53115 Bonn, Germany}

\begin{abstract}
We study dislocation networks in the plane using the vectorial phase-field model introduced by Ortiz and coworkers, 
in the limit of small lattice
spacing. We show that, in a scaling regime where the total length of the dislocations
is large, the phase field model reduces to a simpler model of the strain-gradient type. 
The limiting model contains a term describing the three-dimensional elastic energy
and a strain-gradient term describing the energy of the geometrically necessary dislocations,
characterized by the tangential gradient of the slip. The energy density appearing in the strain-gradient term is determined by the solution of a cell problem, which depends on the line tension energy of dislocations. In the case of cubic crystals with isotropic elasticity { our model shows that}
 complex microstructures may form, in which dislocations with different Burgers vector and orientation
react with each other to reduce the total self energy.
\end{abstract}


\end{frontmatter}

%\linenumbers

\section{Introduction}
Michael Ortiz and coworkers \cite{Ortiz1999,KoslowskiCuitinoOrtiz2002,KoslowskiOrtiz2004}
proposed the use of a vector-valued phase field as a device for describing complex dislocation arrangements.
Their model permits to study situations in which multiple slip systems are active, as long
as the activity is limited to a single slip plane. It incorporates both a local Peierls interplanar
nonconvex potential, which characterizes the discrete nature of slip, and long-range elastic energy. 
Numerical simulations  permitted to identify stable dislocation
structures in finite twist boundaries \cite{KoslowskiOrtiz2004}. The optimal structures obtained from the simulations exhibit
a pattern containing  regular square or hexagonal dislocation networks, 
separated by  complex dislocation pile-ups.

The classical analysis of dislocations is based on regularized continuum models, see \cite{HirthLothe1968,HullBacon}.
The need for a regularization arises from the $1/r$-divergence of the strain close to
the singularity, and is often implemented either by excluding a small volume
around the core or by smoothing the singularity,
in both cases on a length scale of the order of the lattice spacing $\eps$.
In reality, in a discrete lattice there is no singularity, and indeed
the mathematical analysis of dislocation models has shown that the regularization in 
continuum models
plays the same role as the lattice spacing in discrete ones. 
The Ortiz phase-field model, 
as well as the Nabarro-Peierls model, can be
seen as a different way of regularizing linear elasticity.
{The Nabarro-Peierls model is often understood to be a one-dimensional model 
for straight dislocations, but natural extensions to curved dislocations have permitted
to study the energetics of dislocation loops, see for example \cite{XuOrtiz1993,XuArgon2000,XiangWeiMingE2008}.}

The mathematical analysis of the phase-field model 
highlights the occurrence of microstructures over many different length scales. 
Focusing on the regime where the leading-order contribution to the total energy
is given by the dislocation line tension, a number of mathematical papers
rigorously characterized the asymptotics of the model within the framework of 
$\Gamma$-convergence.
This was started in \cite{GarroniMueller2005,GarroniMueller2006}
for the scalar case in which all dislocations have the same Burgers vector, 
then extended in \cite{CacaceGarroni2009,ContiGarroniMueller2011}
to the vectorial situation with multiple slip, and 
in \cite{ContiGladbach2015} to dislocations localized to two parallel planes.
One key result is that
straight dislocations
with certain Burgers vectors and orientations 
may spontaneously decompose into several parallel dislocations, and in some
cases a zig-zag structure is optimal,
see \cite{CacaceGarroni2009,ContiGarroniMassaccesi2015,ContiGladbach2015}.
{ These mathematical results gave a rigorous foundation to the classical Frank's rule for
dislocation reaction \cite{HirthLothe1968,HullBacon}.}

A fully three-dimensional discrete model for dislocations and plasticity was
proposed by Ariza and Ortiz \cite{ArizaOrtiz06}, see also \cite{RamaArizaOrtiz2007}. Their model offers a general 
framework for dislocations in a lattice, and is amenable to a simple analysis
of the continuum limit in situations where Fourier methods are appropriate. 
In the line-tension scaling, more complex techniques are however necessary to
pass to the continuum limit. A rigorous derivation of a line-tension model
from linear elasticity with core regularization was given in \cite{ContiGarroniOrtiz2015},
an extension to a discrete model of the Ariza-Ortiz type will appear in \cite{ContiGarroniOrtizDiscrete}.
Also in this case, relaxation of straight dislocations may be observed, leading to a
line-tension energy which may be smaller \cite{ContiGarroniMassaccesi2015,ContiGarroniOrtiz2015}
than the one predicted by the classical prelogarithmic factor
based on an {\sl ad hoc} generalized plane-strain {\sl ansatz} \cite{BarnettSwanger1971,Rice:1985b}. 

Energy relaxation by formation of microstructure may be even more relevant 
in a situation in which one studies the collective behavior of many dislocations. 
Precisely, one considers a situation in which the total length of the dislocation
lines diverges, and one observes a continuous, macroscopic distribution of dislocations.
Whereas one can estimate the energy of an average dislocation density by adding the
energies of the individual dislocations, interaction and relaxation effects may
significantly alter the picture. 
{ This is a well-known effect in the phenomenological study
of low-angle grain boundaries, see for example \cite{HirthLothe1968,Gottstein2014materialwissenschaft}. % fig. 3.27
We give here a general formulation and a mathematically rigorous 
treatment. In particular, }
we show in Section \ref{sec:cell}
that  { in specific geometries} dislocations with different orientations
and Burgers vectors may interact, {leading in some cases to complex zig-zag patterns}. 
Geometrically, the total (tensorial) density of dislocations corresponds to the 
total incompatibility of the elastic strain field, and therefore to its (distributional) curl. 
For this reason, the energy of  a dislocation density plays a fundamental role in 
the regularization of models of crystal plasticity, leading to the so-called strain-gradient plasticity models
\cite{FleckHutchinson1993,NixGao1998,FleckHutchinson2001,Bassani2001,ContiOrtiz05,KurodaTvergaard2008a,FokouaContiOrtiz2014}.
Indeed, the presence of large latent hardening renders the variational problem of crystal plasticity,
within the deformation theory, nonconvex, leading to lack of existence of minimizers due to the
spontaneous formation of 
fine structure, such as slip bands \cite{AubryOrtiz03,OrtizRepetto1999}, which can again be treated by the theory of relaxation
\cite{ContiOrtiz05,ContiTheil2005,ContiDolzmannKreisbeck2013b,AnguigeDondl2014}.
{The phase-field model by Ortiz and coworkers that we study here was related to strain gradient plasticity
in \cite{HunterKoslowski2008}, where in particular the parameters of continuum strain-gradient plasticity approach
were derived from the phase-field dislocation model.}
{The corresponding process for the forces is the derivation of a  continuum approximation 
for the Peach-Kohler force, as derived in \cite{Xiang2009,ZhuXiang2014}. We remark that in all these works
the {\em relaxation} of the dislocation structures, which naturally arises if a  mathematically rigorous
variational limiting procedure is attempted, is not considered.}

Strain-gradient plasticity models can be rigorously derived from discrete models,
or regularized semidiscrete models, using $\Gamma$-convergence with a choice of the scaling of the energy
which balances the contributions of the elastic field and of the dislocation core energies. 
This was performed for the first time  for point
dislocations in the plane by Garroni, Leoni and Ponsiglione
\cite{GarroniLeoniPonsiglione2010}
 in a geometrically linear setting with a core regularization approach, 
 and by M\"uller, Scardia and Zeppieri with a geometrically nonlinear formulation
\cite{MuellerScardiaZeppieri,MuellerScardiaZeppieri2015}. Both results rely on a well-separation assumption, which permits to locally
estimate the self-energy of each individual dislocation. The assumption of  point dislocations in the plane
corresponds to an array
of straight parallel dislocations in three dimensions.

In this work, we derive a strain-gradient model for a density of line dislocations
in the plane. Our starting point is the vectorial phase-field model developed by Ortiz and coworkers.
Our key result is that the energy can be approximated by the sum of two terms, given by
the long-range elastic interaction and the self-energy of the dislocation density,
see (\ref{eqdeffzero}) below. The self-energy
itself is determined by solving a cell problem, 
{ which corresponds to selecting the energy of the optimal dislocation structure among all those
with the same average dislocation density. In particular, it is in general smaller than the sum of the 
line-tension energies obtained using individual straight dislocations.}
We remark that the key ingredients in this relaxation is the anisotropy of the 
prelogarithmic factor in the energy of a single, straight dislocation. Higher-order interaction 
effects may further enrich the picture.

We introduce in Section \ref{sec:model} the vectorial phase-field model on a torus that we shall use for the rest of the paper
and the relevant scaling regime.
The line-tension energy of individual dislocations is discussed in Section \ref{sec:psi}, and the energy
of dislocation structures in Section \ref{sec:cell}.
Finally, in Section \ref{sec:result} we present the full limiting model which contains both elasticity and dislocation self-energy.

\section{The vectorial phase-field model}\label{sec:model}
Following  Ortiz and coworkers \cite{Ortiz1999,KoslowskiCuitinoOrtiz2002,KoslowskiOrtiz2004}, we study dislocation
patterns in the plane, assuming periodicity in the transversal directions. 
For $L>0$ we consider the domain $\T=(0,L)^2$ with periodic boundary conditions. 
The elastic deformation $U:\T\times \R\to\R^3$ is equally assumed to be periodic in the horizontal variables, and may jump across the $\{x_3=0\}$ plane.
The elastic energy of $U$ is  complemented by an additional term due to the short-range interatomic interactions across 
the slip plane, resulting in the total energy
\begin{equation}\label{eqdefhate}
F_\eps[U]=\frac1\eps\int_\T W_\calB(\gamma) dx+ \int_{\T\times \R} \frac12 \C e(U)\cdot e(U) \, dx\,.
\end{equation}
Here $e(U)=(\nabla U+\nabla U^T)/2$ is the elastic strain, $\calB\subset\R^2\times\{0\}$ is the two-dimensional lattice of possible slip
vectors, $U:\T\times\R\to\R^3$ is the $(0,L)^2$-periodic displacement field, and $\gamma=[U]$ is its jump across the $\{x_3=0\}$ plane, i.e., the plastic slip,
{often called  disregistry in the context of Peierls-Nabarro models.}
The latter is assumed to take values in the linear space generated by $\calB$, which typically is $\R^2\times\{0\}$, reflecting volume conservation.
The potential $W_\calB:\R^2\times\{0\}\to[0,\infty)$ vanishes on all vectors $v\in\calB$ and is positive elsewhere.
Further, $\C:\R^{3\times 3}\to \R^{3\times 3}$ is the (symmetric) tensor of linear elastic coefficients, which obeys
for some $c>0$ the conditions
\begin{equation}\label{eqasscelas}
 \frac1c|\xi+\xi^T|^2\le \C \xi\cdot \xi \le c |\xi+\xi^T|^2 \text{ for all } \xi\in \R^{3\times 3}\,.
\end{equation}
The presence of the large coefficient $1/\eps$ in the first term in (\ref{eqdefhate}) 
can be understood as a remnant of the fact that in a discrete model the first term accounts
only for interactions across a plane, and that $U$ should be understood as a displacement divided by the 
lattice parameter $\eps$. We refer to \cite{KoslowskiOrtiz2004,GarroniMueller2005} 
for a more detailed discussion of this model.

In the following we are specifically interested in the slip field $\gamma=[U]:\T\to\R^2\times\{0\}$. 
Let $s_1,\dots, s_N\subset\calB$ be a set of Burgers vectors which forms a basis for $\calB$. 
Then one can express $[U]$ as a linear combination
of $s_1,\dots, s_N$, with coefficients given by a map $u:\T\to\R^N$,
\begin{equation}\label{eqgammau}
\gamma(x)=[U](x)=\sum_{i=1}^N u_i(x)s_i \,.
 \end{equation}
Minimizing out
the displacement field $U$ for fixed $u$ with the aid of Fourier transform,
 and approximating $W_\calB$ by a piecewise quadratic Peierls potential,
leads to
\begin{equation}
\begin{split}
 E_\eps[u] =&\frac1\eps\int_\T \dist^2(u,\Z^N) dx\\
 &+ \sum_{i,j=1}^N\int_{\T\times \T} K_{ij}(z)(u_j(x)-u_j(x+z))(u_i(x)-u_i(x+z)) \, dx\, dz\,.
 \end{split}
\end{equation}
Although physically $N=2$ is the most relevant case, we keep the
dimension $N$ general, for easier comparison with the scalar $N=1$ case.
The interaction kernel $K:\T\to\R^{N\times N}$ is symmetric, in the sense that {$K^T(\xi)=K(\xi)$}, and homogeneous of degree $1$ in Fourier space, 
in the sense that its Fourier coefficients obey $\hat K(\xi)=|\xi|\hat K(\xi/|\xi|)$. Transforming back to real space shows that
\begin{equation}\label{eqkgammar}
 K(z)=\Gamma(z)+R(z)\,,
\end{equation}
where $\Gamma$ denotes the singular part, which is homogeneous of degree $-3$, 
\begin{equation}
\label{eq:derGammaintro}
\Gamma(z)=\frac{1}{|z|^{3}} \hat\Gamma\left(\frac{z}{|z|}\right)\,,
\end{equation}
where $\hat\Gamma:S^1\to\R^{N\times N}$ 
{is assumed to obey, for some $c>0$, }
%(related to, but not identical with, the constant $c$ in (\ref{eqasscelas}))
\begin{equation}\label{eqgammaposdeftheointro}
\frac{1}{c}|v|^2\le v\cdot \hat\Gamma(z)v\le c|v|^2 \hskip5mm 
\text{ for all } v\in \R^N, \, z\in S^1\,.
\end{equation}
The specific form depends on the elastic constants of the crystal, for example in an elastically
isotropic cubic crystal one has 
\begin{equation}\label{eqgammacubic}
 \Gamma^\cubic(z)=\frac{\mu}{16\pi(1-\nu) |z|^3}\left(\begin{matrix}\nu+1-3\nu\frac{z_2^2}{|z|^2} &   3\nu\frac{z_1z_2}{|z|^2}\\  3\nu\frac{z_1z_2}{|z|^2}&\nu+1-3\nu\frac{z_1^2}{|z|^2}  \end{matrix} \right)\,,
\end{equation}
where $\nu$ and $\mu$ denote the 
material's Poisson's ratio and shear modulus, respectively {(see \cite{CacaceGarroni2009}).}
{
It is easy to see that for $\mu>0$ and $\nu\in (-1,1/2)$ the kernel $ \Gamma^\cubic$ fulfills the assumption (\ref{eqgammaposdeftheointro}).}

The correction $R$ in (\ref{eqkgammar}) 
is called the regular part of the interaction, which arises from the periodic boundary conditions and it is bounded,
in the sense that there is $c>0$ such that $|R(z)|\le c$ for all $z$. The total energy density
is nonnegative for every $z$, in the sense that $K(z)v\cdot v\ge 0$ for all $z\in \T$, $v\in \R^N$.

The singular part of the energy will be particularly important for our analysis.
We denote it by 
\begin{equation}
 \calE[u]=
 \sum_{i,j=1}^N\int_{\T\times \T} \Gamma_{ij}(z)(u_j(x)-u_j(x+z))(u_i(x)-u_i(x+z)) \, dx\, dz\,.
\end{equation}
\begin{figure}
 \begin{center}
  \includegraphics[width=9cm]{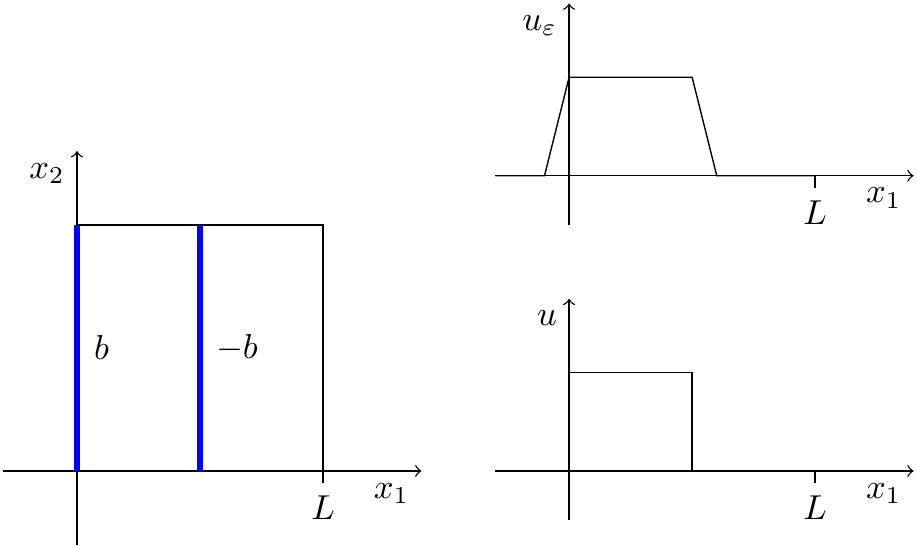}
 \end{center}
\caption{Sketch of the simple situation with two dislocations with {Burgers vector} $b$ and $-b$ 
along $\R e_2$ and $\R e_2+L/2 e_1$. The right panel illustrates the corresponding 
values of $u$ and $u_\eps$, as defined in (\ref{eqdefu}) and (\ref{eqdefueps}).
In all cases, only one period is shown.
}
\label{fig-dislo}
\end{figure}
In order to understand the appropriate energy scaling, we start from a simple one-dimensional situation. Assume that a 
dislocation with Burgers vector $b$ is given, along the line $\R e_2$ (which is the $x_2$ axis), and 
one with Burgers vector $-b$ along the parallel line $\R e_2+L/2 e_1$, see Figure \ref{fig-dislo}. 
The sharp-interface model would then have
\begin{equation}\label{eqdefu}
 u(x)=\begin{cases}
       1 & \text{ if } kL<x_1<(k+\frac12) L\text{ for some } k\in\Z\,,\\
       0 & \text{ otherwise}.
      \end{cases}
\end{equation}
It is easy to see that the energy of this slip is infinite, as the elastic energy diverges. A natural regularization on
the scale $\eps$ is 
\begin{equation}\label{eqdefueps}
 u_\eps(x)=
 \begin{cases}
 \max \{0, 1- \frac1\eps \dist(x_1, [0,L/2])\} & \text{ for } 0\le x_1<L\\
 \text{periodic extension } & \text{ otherwise.}
      \end{cases}
\end{equation}
One can then easily compute that $E_\eps[u_\eps]\sim \ln 1/\eps$. This corresponds to the well-known logarithmic 
divergence of the energy per unit length of dislocations.

We consider now a situation in which $M$ such dislocations are present. The total energy will be of order
$M\ln \frac1\eps$, the total variation of $u_\eps$ will be of order $M$, since it has $M$ ``jumps'' of height 1.
Correspondingly, the elastic energy behaves as $M^2$. Therefore if $M\sim \ln \frac1\eps$ the total line-tension energy and the
total elastic energy are of the same order, see \cite{GarroniLeoniPonsiglione2010,MuellerScardiaZeppieri,MuellerScardiaZeppieri2015} for mathematically
rigorous treatments of this heuristics.

For this reason, in the following we shall focus on a situation in which 
$E_\eps[u_\eps]$ is of order $(\ln \frac1\eps)^2$, and $u_\eps$ is itself of order 
$\ln \frac1\eps$.  Specifically, we are interested in the limit $\eps\to0$, assuming that
\begin{equation}\label{eqconvu}
%\lim_{\eps\to0}
\frac{u_\eps}{\ln \frac1\eps}\to u_0
\end{equation}
and computing  the asymptotic energy (in the sense of $\Gamma$-convergence, see below)
\begin{equation}\label{eqconveeps}
%\lim_{\eps\to0}
\frac{E_\eps[u_\eps]}{(\ln \frac1\eps)^2}\to E_0[u_0]\,.
\end{equation}
The limiting energy $E_0$ will turn out to contain both a long-range elastic energy term and 
a short-range self-energy term, which characterizes the planar distribution of dislocations. 
The limiting procedure should be understood as determining the effective energy $E_0$ of a limiting
plastic slip distribution $u_0$ by approximating it with an optimal sequence $u_\eps$, and the associated distribution of dislocations, and computing the optimal
energy along the sequence.

\section{Dislocation line-tension energy}\label{sec:psi}
The starting point of our analysis is the line-tension energy approximation derived in 
\cite{CacaceGarroni2009,ContiGarroniMueller2011}. We now briefly review some results
which will be needed in the following.

The prelogarithmic factor $\psi_0(b,n)$ of a dislocation with Burgers vector $b\in\Z^N$ 
(in the coordinates given in (\ref{eqgammau}) above) can be computed as
\begin{equation}\label{eqdefpsi0}
 \psi_0(b,n)=2\int_\R \Gamma(n+tn^\perp)b\cdot b\, dt
\end{equation}
where $n^\perp=(-n_2,n_1)$ is the orthogonal vector to $n$. The vector $n$, in turn, is a unit vector normal to the dislocation line,
so that $n^\perp$ is the tangential vector.  For an elastically isotropic crystal with cubic symmetry, one obtains 
\begin{equation}\label{eqdefpsi0cubic}
 \psi_0^\cubic(b,n)=\frac{\mu }{4\pi} \left(|b|^2+\eta (b\cdot n)^2\right)\,,
\end{equation}
where $\eta=\nu/(1-\nu)$ is a material parameter which depends on the material's Poisson's ratio $\nu$
(see \cite[Eq. (51)]{KoslowskiCuitinoOrtiz2002} or 
{\cite[Eq. (4.2)]{ContiGarroniMassaccesi2015}, the latter is missing a factor $1/2$; 
in  \cite{CacaceGarroni2009,ContiGarroniMueller2011} a term $\mu/2$ was factored out}).
The parameter $\eta$ characterizes the relative energy difference between a screw
dislocation, which has $b$ orthogonal to $n$, and an edge dislocation, which has $b$ parallel to $n$.
{The expression (\ref{eqdefpsi0cubic}) agrees with the classical energy per unit length
of dislocations, as given for example in \cite[Eq. (3.13) and (3.52)]{HirthLothe1968}.}

By prelogarithmic factor we mean that the energy per unit length of the dislocation is, to leading order,
$\ln \frac1\eps\, \psi_0(b,n)$. In the simple case $b=e_1$ and $n=e_1$,
the slip field defined in (\ref{eqdefueps}) has indeed {to leading order} energy 
$\ln \frac1\eps\, \psi_0(e_1,e_1)$.
The reduction of the integration from two to one dimension is made using (\ref{eq:derGammaintro}), 
see \cite{GarroniMueller2005,ContiGarroniMueller2011} for details.
The prelogarithmic factor $\psi_0$ can also be computed directly starting from three-dimensional elasticity,
using for example a core-radius regularization. The explicit computation can be done with a suitable
plane strain ansatz \cite{HirthLothe1968}; a variational characterization minimizing the energy
in long cylinders can be found in \cite{ContiGarroniOrtiz2015}. All these methods give the same result for $\psi_0$,
which only depends on the matrix of elastic constants $\C$.

\begin{figure}
 \begin{center}
  \includegraphics[width=12cm]{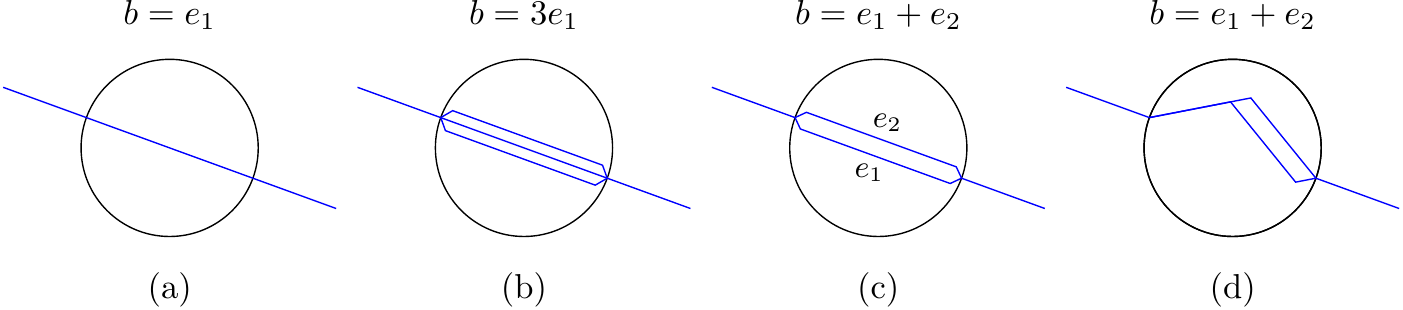}
 \end{center}
\caption{Sketch of {some} possible microstructures arising in the definition of $\psi_\rel$ (\ref{eqdefpsirel}).
These microstructure arise at the intermediate scale between the one of the lattice and the macroscopic one, 
as in Figure \ref{fig-multiscale}(b).}
\label{fig7}
\end{figure}

The prelogarithmic factor $\psi_0$ gives, after rescaling, an effective energy per unit length. One can therefore
define a line-tension model for a dislocation network. A dislocation network can be parametrized by finitely many
oriented curves $\gamma_1,\dots,\gamma_M\subset\T$, 
each with an associated Burgers vector $b_1,\dots, b_M$, which is a conserved
quantity in the sense that for every point $x$, where one or more curves start or end, the sum of the Burgers vectors
of the curves ending at $x$ equals the sum of the Burgers vectors of the curves starting at $x$. The (unrelaxed) 
line-tension energy is then
\begin{equation}\label{eqenergypsi0}
 I_0=\sum_{i=1}^M \int_{\gamma_i} \psi_0(b_i, n) ds
\end{equation}
where $ds$ denotes integration along the curve and $n$ is the normal vector to the curve.

\begin{figure}
 \begin{center}
  \includegraphics[width=9cm]{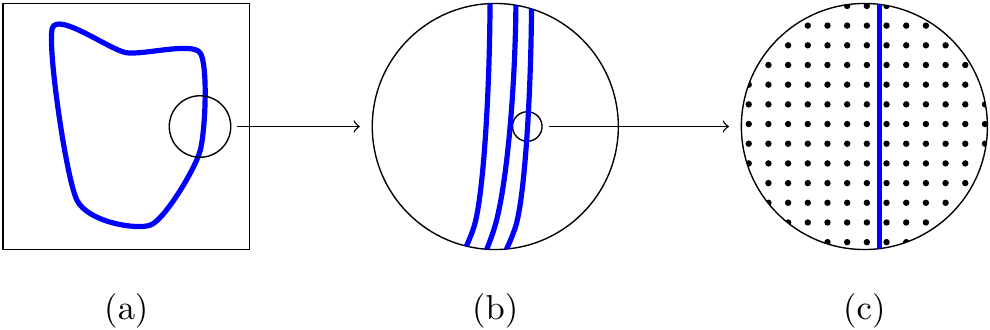}
 \end{center}
\caption{Multiscale relaxation of the dislocation energy. (a) shows a macroscopic dislocation line, whose energy
per unit length is characterized by $\psi_\rel$. (b) illustrates a blow-up of a small portion of that dislocation,
which - in this example - is subdivided into three separate dislocation lines. These lines are very close on the scale of the sample, so that they are not distinguished in (a), but well separated on the scale of the lattice, so that the total energy is the  sum $\psi_0(b_1,n)+\psi_0(b_2,n)+\psi_0(b_3,n)$.
(c) illustrates a further blow-up, on the scale of the lattice, where only one of the three dislocation
lines from (b) is seen.
}
\label{fig-multiscale}
\end{figure}

The energy density $\psi_0(b,n)$ is quadratic in $b$, as is apparent 
 from the expression in (\ref{eqdefpsi0}). 
 This does not, however, reflect the true macroscopic energetic cost of a singularity with { total} Burgers vector $b$.
Consider for example the case that the curve $\gamma_1$ has { total} Burgers vector $b_1=(3,0)$. The
energy of this is 9 times the energy of a dislocation with the same path and Burgers vector $(1,0)$,
since $\psi_0(3e_1,n)=9\psi_0(e_1,n)$. 
It is therefore energetically convenient to split $\gamma_1$ into three dislocations with smaller Burgers vectors,
say $\gamma_1'$, $\gamma_1''$, $\gamma_1'''$, with $b_1'=b_1''=b_1'''=(1,0)$. The three curves have the same
start and endpoint as $\gamma_1$, but are otherwise disjoint, although very close to each other, and have,
up to higher order terms, the same length as $\gamma_1$. The total energy along these three curves 
is then only three times
the energy of the original curve (see Figure \ref{fig7}b)
This mathematical observation corresponds to the physically well-known fact that only the shortest lattice vectors are stable Burgers vectors of dislocations.
{ In general, this shows that an effective dislocation energy can be stable with respect to decay into parallel dislocations 
only if it is subattidive in the first argument, a condition which corresponds to the classical 
Frank's rule for dislocation reactions \cite{HirthLothe1968,HullBacon}.}
It is important to notice that when we deal with the line tension energy we are considering dislocations at a scale much larger than the lattice spacing, where the core region is identified with a line. Therefore a configuration with a larger { total}  Burgers vector, as the one considered in this examples, that might look unphysical, needs to be understood as a cluster of parallel dislocation lines with short Burgers vectors.
The energy $\psi_0(3e_1,n)$ describes the (hypothetical) situation in which the individual lines have a 
separation of the order of the lattice spacing, so that the total elastic energy of the three dislocations  
is $9\psi_0(e_1,n)$. The energy $3\psi_0(e_1,n)$ can be achieved if the three dislocations
are so close that they can still be identified at a macroscopic scale, but their 
relative distance is large enough to avoid interaction at the leading order, resulting on an elastic distortion which is the superposition of the effect due to each single dislocation.
{If the energy computation were attempted numerically, it would be important that the resolution is high enough to resolve the separation of the curves, i.e., the microstructure. In doing the relaxation step
analytically, the infinite-resolution limit is automatically incorporated.} 
See Figure \ref{fig-multiscale} for an illustration of the different scales.

A less evident instance of the same  effect arises for linear combinations of different Burgers vectors.
For example, a curve with $b_1=(1,1)$ can be replaced by a curve with $b_1'=(1,0)$ and another one with $b_1''=(0,1)$,
 see Figure \ref{fig7}c.
Even more, a straight curve may be replaced by a finely-oscillating curve, which has more length but possibly
an energetically  more convenient orientation (like in faceting of crystal surfaces), see Figure \ref{fig7}d.
{ Whether this is energetically convenient depends on the details of the problem, including in particular the orientation of the dislocation line and the material's elastic constants.}
For a more detailed discussion of these phenomena we refer to \cite{CacaceGarroni2009,ContiGarroniOrtiz2015},
and to \cite{AmbrosioBraides1990a,AmbrosioBraides1990b,ContiGarroniMassaccesi2015} for the mathematical background.

The optimal energy per unit length of a dislocation network that carries the { total}  Burgers vector $b\in\Z^N$ across
 a segment with normal $n\in S^1$ is given by the cell-problem formula
 \begin{equation}
  \psi_\rel(b,n)=\inf \Bigl\{\sum_i \int_{\gamma_i} \psi_0(b_i, \nu) ds \Bigr\}
 \end{equation}
where the minimum is taken over all networks as described above, which start in the point 0 with a Burgers vector $b$, 
and end in the point $n^\perp$ with the same Burgers vector, {$\nu$ is the normal to the dislocation line $\gamma_i$}.
This minimization corresponds to the optimization among all possible dislocation structures which are admissible, in the sense
that the { total}  Burgers vector is conserved, and which agree with a straight dislocation with { total}  Burgers vector $b$ and normal $n$ outside a small region, as illustrated in Figure \ref{fig7}.

An equivalent formulation can be given in terms of piecewise constant phase fields, which correspond to 
functions of bounded variation with values in $\Z^N$. Indeed,
 \begin{equation}\label{eqdefpsirel}
  \psi_\rel(b,n)=\inf \Bigl\{\frac12\int_{J_u}  \psi_0([u], \nu_u) d\calH^1: u\in BV(B_{1},\Z^N), \, u=u^{b,n} \text{ on } 
  \partial B_{1}\Bigr\}
 \end{equation}
where $u^{b,n}(x)=0$ if $x\cdot n<0$, and $b$ if $x\cdot n>0$, and 
 $B_{1}$ the ball of unit radius centered in the origin, which has diameter $2$.
Here $J_u$ is the set where $u$ jumps, which corresponds to the union of the dislocation curves $\gamma_i$, $[u]$ is the jump,
which corresponds to the local Burgers vector, {and $\nu_u$ the normal to $J_u$. }
For $u^{b,n}$, the jump set is the diameter of $B_{1}$ orthogonal to $n$, the jump equals $b$, and therefore
setting $u=u^{b,n}$ in (\ref{eqdefpsirel}) yields $\psi_\rel(b,n)\le \psi_0(b,n)$.
Here and below $\calH^1$ denotes integration along the (one-dimensional) set $J_u$.
We refer to \cite{AmbrosioFP} for precise definitions
of these concepts. 

By \cite[Lemma 6.4]{ContiGarroniOrtiz2015} one can show that 
 \begin{equation}\label{eqdefpsirelpc}
  \psi_\rel(b,n)=\inf\Bigl\{\frac12\int_{J_u}  \psi_0([u], \nu_u) d\calH^1: u\in PC(B_{1},\Z^N), \, u=u^{b,n} \text{ on } 
   \partial B_{1}\Bigr\}
 \end{equation}
where, given $\omega\subset\R^2$ and a set $A$, we denote by $PC(\omega;A)$ the set of 
polygonal piecewise affine functions with values in $A$, i.e., the set of functions $u:\omega\to A$ such that $\omega$
can be covered by finitely many disjoint polygons, such that $u$ is constant on each of them.
In particular, the construction in Figure \ref{fig7}(b) permits to prove that $\psi_\rel(b,n)$ has linear growth in $b$, in the sense that
\begin{equation}\label{eqpsirellingrw}
 \frac1c |b| \le \psi_\rel(b,n) \le c|b|\text{ for all }n\in S^1,b\in\Z^N\,.
\end{equation}

Starting from Michael Ortiz's phase field model, 
in \cite{GarroniMueller2005,CacaceGarroni2009,ContiGarroniMueller2011} it was shown that,
in the regime in which the energy is proportional to $\ln \frac1\eps$, the energetically optimal 
slips correspond to a dislocation distribution with finite total length, whose energy can be computed at 
leading order using $\psi_\rel$. Precisely, this can be expressed in terms of $\Gamma$-convergence as
\begin{equation}\label{eqgammaconveepse0}
 \frac{1}{\ln (1/\eps)} E_\eps \stackrel{\Gamma}{\to} E_0^*
\end{equation}
where 
\begin{equation}
 E_0^*[u]=\int_{J_u} \psi_\rel([u],n) d\calH^1 \text{ for } u\in BV(\T;\Z^N)\,.
\end{equation}
The convergence in (\ref{eqgammaconveepse0}) means that for any sequence $u_\eps\to u_0$
one has
\begin{equation}
 E_0^*[u_0]\le \liminf_{\eps\to0} \frac{1}{\ln (1/\eps)} E_\eps[u_\eps] 
\end{equation}
and, conversely, that for any $u_0$ there is a sequence 
$u_\eps\to u_0$
such that
\begin{equation}
 E_0^*[u_0]= \lim_{\eps\to0} \frac{1}{\ln (1/\eps)} E_\eps[u_\eps] \,,
\end{equation}
see \cite{Dalmaso1993,Braides02}.
This corresponds to the fact that $E_0^*[u_0]$ is the smallest energy, to leading order,
among all sequences $u_\eps$ converging to $u_0$ in $L^1(\T;\R^N)$. 

The scaling 
of the energy in (\ref{eqgammaconveepse0}) 
and the convergence $u_\eps\to u_0$  used in 
\cite{GarroniMueller2005,CacaceGarroni2009,ContiGarroniMueller2011}  are appropriate to emphasize
the line-tension energy of individual dislocation lines. In contrast, in the present work we focus
on a different regime, characterized by (\ref{eqconvu}) and (\ref{eqconveeps}),  which  is appropriate for
computing the total energy of a complex dislocation structure, consisting of a large number of individual dislocation lines,
as will be explained in the next section.

\section{Cell structures and their energy}
\label{sec:cell}

\begin{figure}
 \begin{center}
  \includegraphics[width=9cm]{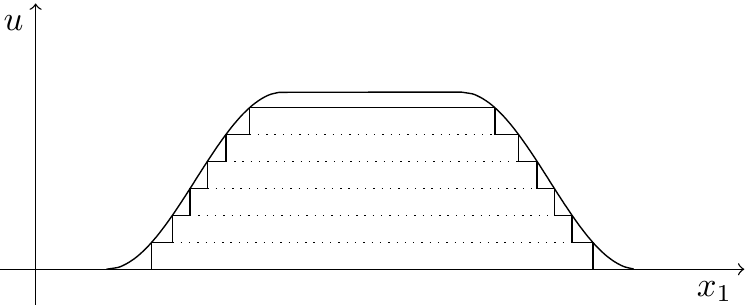}
 \end{center}
\caption{Approximation of a continuous slip by a step process with many individual dislocations.}
\label{fig:scala}
\end{figure}

\begin{figure}
 \begin{center}
  \includegraphics[width=6cm]{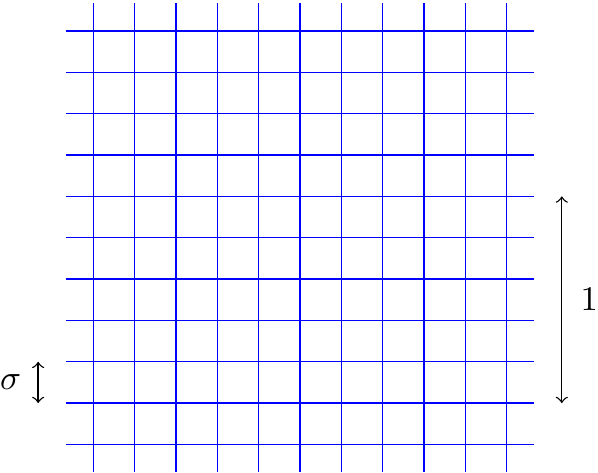}
 \end{center}
\caption{Sketch of possible simple grid structures as defined in (\ref{eqdefhatu}) and (\ref{eqdefusigma}),
corresponding to the macroscopic slip $u(x)=Ax=(x_1,-x_2)$.}
\label{fig:griglia}
\end{figure}

%\subsection{Cell structures and their energy}
We now study macroscopically continuous slip. Passing to a larger scale, one sees a large number of
dislocations, with a large total Burgers vector, which may be approximately uniformly distributed in space, see Figure \ref{fig:scala}.
{This is the same procedure that is usually used to study low-angle grain boundaries or the opening of cracks, see
for example \cite{XuOrtiz1993,BlatovKaxiras1997,DaiXiangSrolovitz2013}}.
%After rescaling, this can be understood as the superposition of many dislocations with small Burgers vector.
Consider for definiteness a locally affine slip, $v(x)=Ax$, for some matrix $A\in \R^{N\times 2}$.
Before presenting the general situation, we start from the illustrative example $A=e_1\otimes e_1-e_2\otimes e_2$.
%Denoting by $\sigma$ the rescaling parameter, this needs to be approximated with dislocations with
%Burgers vector in $\sigma\Z^N$.
We write
\begin{equation}\label{eqdefhatu}
\hat u_\sigma(x)= e_1  \left\lfloor \frac{x_1}{\sigma}\right\rfloor -
e_2 \left\lfloor \frac{x_2}{\sigma}\right\rfloor\,,
\end{equation}
where $\sigma$ is a small parameter and
$\lfloor y\rfloor=\max\{z\in\Z: z\le y\}$ denotes the largest integer not larger than $y$.
This is a piecewise constant function, which for small $\sigma$ is close to $A x/\sigma$
and jumps across horizontal and vertical lines spaced by $\sigma$, which represent the dislocations. The amplitudes of the jumps are
$e_2$ and $e_1$, respectively (see Figure \ref{fig:griglia}).
The energy per unit area of this configuration can be computed using (\ref{eqenergypsi0}) and is $\frac1\sigma\psi_0(e_1,e_1)+
\frac1\sigma\psi_0(e_2,e_2)$. As above, the energy can in principle be reduced by local microstructures, leading to relaxation, the corresponding energy will be obtained using $\psi_\rel$ instead of $\psi_0$.

Since both the slip and the energy are diverging as $\sigma$ tends to zero, it is convenient to rescale.
We multiply both the slip and the energy by $\sigma$ and obtain the sequence of slip fields
\begin{equation}\label{eqdefusigma}
u_\sigma(x)= e_1\sigma  \left\lfloor \frac{x_1}{\sigma}\right\rfloor -
e_2 \sigma\left\lfloor \frac{x_2}{\sigma}\right\rfloor
\end{equation}
converging to $Ax$ for small $\sigma$. Its jump set has (locally) length diverging as $1/\sigma$, but its 
jumps have amplitude $\sigma$. It is then natural to correspondingly rescale the energy density, and to consider
\begin{equation}\label{eqdefpsiinfty}
{\psi_\infty(b,n)=
  \liminf_{\sigma\to0} \sigma \,\psi_\rel(\frac b\sigma,n) 
=
  \liminf_{s\to\infty} \frac{1}{s} \psi_\rel(s b,n) \,.}
\end{equation}
This limit is finite since $\psi_\rel$ has linear growth (recall (\ref{eqpsirellingrw})), the new
function $\psi_\infty$ captures the asymptotic behavior of $\psi_\rel$ at infinity, i.e., the approximate energy per unit
length and unit Burgers vector. 

The map $u_\sigma$ jumps, within the unit square $[0,1)^2$, on $\lceil 1/\sigma\rceil$ horizontal segments of
unit length, 
and the jump amplitude is $\sigma e_2$
(with $\lceil y\rceil=\min\{z\in\Z: z\ge y\}$, {since $\psi_0^\cubic$ is an even function the orientation does not matter}).
Their (scaled) total energy is $\lceil\frac1\sigma\rceil \psi_\infty(\sigma e_2,e_2)$. {The same holds for the jumps over horizontal lines}, and
 the energy per unit area of $u_\sigma$ is
 \begin{equation}
 \lceil\frac1\sigma\rceil \psi_\infty(\sigma e_1,e_1)+
\lceil\frac1\sigma \rceil\psi_\infty(\sigma e_2,e_2)\sim \psi_\rel(e_1,e_1)
+\psi_\rel(e_2,e_2)\,.
\end{equation}

\begin{figure}
 \begin{center}
  \includegraphics[width=12cm]{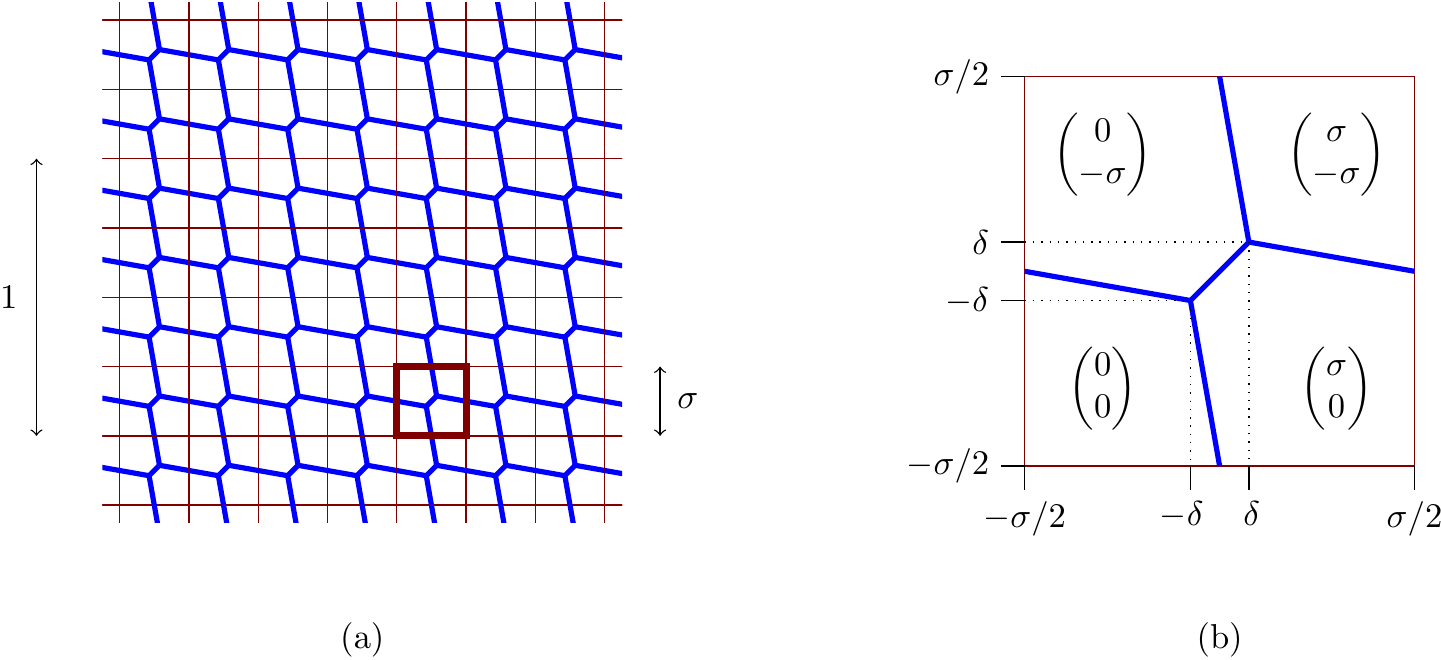}
 \end{center}
\caption{Sketch of a more complex microstructure
with the energy given in 
(\ref{eqdefedelta}).
(a) shows a larger region, (b) a blow-up of a $\sigma\times\sigma$ square on which the computation 
is performed.
For $\delta=0$ this reduces to the microstructure in Figure \ref{fig:griglia}.}
\label{fig-micro1}
\end{figure}

In many cases, however, more complicated structures appear. Consider for example the pattern illustrated in Figure 
\ref{fig-micro1}. Here the two sets of dislocations interact, and overlap over part of the domain 
to form a composite dislocation with a larger Burgers vector parallel to $(1,1)$.
In order to compute the energy of this configuration, we assume that the central segment is
oriented at $45$ degrees and denote by $2\delta\in[0,\sigma)$ the length of its horizontal projection, see Figure \ref{fig-micro1}(b).
The energy of this central segment is $2\delta\sqrt2  \psi_\infty(\sigma(e_1+e_2),{(e_1-e_2)/\sqrt2)}$.
The lower segment has length $\sqrt{\delta^2+(\sigma/2-\delta)^2}$,
normal $\bar n=(\sigma/2-\delta,\delta)/\sqrt{\delta^2+(\sigma/2-\delta)^2}$, 
and Burgers vector $\sigma e_1$. Therefore its energy is
\begin{equation}
 \sqrt{\delta^2+(\sigma/2-\delta)^2}\psi_\infty(\sigma e_1, \bar n)\,.
\end{equation}
Correspondingly, the  segment on the left has energy
\begin{equation}
 \sqrt{\delta^2+(\sigma/2-\delta)^2}\psi_\infty(\sigma e_2,\bar n')
\end{equation}
where $\bar n'= (\delta,\sigma/2-\delta)/\sqrt{\delta^2+(\sigma/2-\delta)^2}$.
The other two are the same, up a translation. 
The energy in the $(0,\sigma)^2$ square is therefore
\begin{align}\nonumber
 e(\delta)=& 2\delta\sqrt2  {\psi_\infty(\sigma(e_1+e_2),(e_1-e_2)/\sqrt2)}\\
 &+\label{eqdefedelta}
 2\sqrt{\delta^2+(\sigma/2-\delta)^2} \left(\psi_\infty(\sigma e_1,\bar n)+\psi_\infty(\sigma e_2,\bar n')\right)\,.
\end{align}
We now show that, at least in the case of a cubic crystal with isotropic elasticity, this configuration is energetically
more convenient than the simple one described in (\ref{eqdefusigma}) and Figure \ref{fig:griglia} above, and which corresponds
to the $\delta=0$ case of the present one.
For cubic crystals we know that $\psi_0$ is given by (\ref{eqdefpsi0cubic}). Further, 
it was shown in \cite[Lemma 4.4]{ContiGarroniMassaccesi2015} that
\begin{equation}
 \psi_\rel^\cubic(ke_i, n)={|k|}\psi^\cubic_0(e_i,n)\text{ for $i=1,2$,  $k\in \Z$ and $n\in S^1$. }
\end{equation}
{From \cite[Lemma 4.6]{ContiGarroniMassaccesi2015} 
using $\psi^\cubic_0(e_1+e_2,(e_1-e_2)/\sqrt2)=2\le \psi^\cubic_0(e_1+e_2,n)$ one easily sees
that}
\begin{equation}
 \psi_\rel^\cubic(k(e_1+e_2), \frac{e_1-e_2}{\sqrt2})={|k|}\psi^\cubic_0((e_1+e_2), \frac{e_1-e_2}{\sqrt2})
 \text{ for  $k\in \Z$  }\,.
\end{equation}
In particular, recalling (\ref{eqdefpsi0cubic}) and (\ref{eqdefpsiinfty}), we obtain
\begin{equation}
 \psi_\infty^\cubic(ke_i,n)=\frac\mu{4\pi} |k| (1+\eta (n\cdot e_i)^2)
\end{equation}
and
\begin{equation}
 \psi_\infty^\cubic(k(e_1+e_2),\frac{e_1-e_2}{\sqrt2})=\frac\mu{4\pi} 2|k|\,.
\end{equation}
Inserting  in (\ref{eqdefedelta})
we obtain
\def\sdelta{\delta}
\begin{align}
e(\delta)&= 4\sigma
\frac\mu{4\pi}
\left[ \sqrt{\sdelta^2+(\frac12\sigma-\sdelta)^2} +\eta \frac{(\frac12\sigma-\sdelta)^2}{\sqrt{\sdelta^2+(\frac12\sigma-\sdelta)^2}} + \sdelta\sqrt2\right]
\\
 &= 4\sigma^2
\frac\mu{4\pi} \left[ \frac{1+\eta}{2} + \frac{\sdelta}{\sigma}(\sqrt2-1-\eta) + O\left(\frac{\sdelta^2}{\sigma^2}\right)\right]\,.
\end{align}
Therefore for $\eta>\sqrt2-1$ the minimum of $e$ is not taken at $\delta=0$. For example, if
 $\eta=1/2$ the optimal value is taken at $\delta\sim 0.07 \sigma$, leading to the pattern illustrated in Figure \ref{fig-micro1}.
Recalling that $\eta=\nu/(1-\nu)$, we can easily see that the condition 
$\eta>\sqrt2-1$ corresponds to
\begin{equation}
 \nu>1-\frac1{\sqrt2}\sim 0.29\,,
\end{equation}
which is satisfied in many metals.
The specific value $\eta=1/2$ chosen for the computation above corresponds to $\nu=1/3$.
{ We remark that this result, although mathematically rigorous under the stated assumptions, 
depends on the chosen orientation and the chosen elastic properties. Further, its physical relevance is restricted
by the assumptions of the current modeling, which in particular focuses on the leading-order energy terms with respect to the small parameter $\eps$.}

%\subsection{Definition of the energy densities}
As illustrated in the above examples, we define the energy per unit area of a general affine slip field
$u(x)=Ax$ by considering the optimal energy of  dislocation networks which realize it.
We start from the rescaled energy density defined in (\ref{eqdefpsiinfty}).
We define $g:\R^{M\times 2}\to[0,\infty)$ by the cell problem
\begin{align}\label{eqdefg}
 g(A)=\inf \Bigl\{ \liminf_{j\to\infty} \frac{1}{\pi}\int_{J_{u_j}\cap  B_1} \psi_\infty([u_j], {\nu_{u_j}}) \,d\calH^1: \ &u_j\text{ piecewise constant,} \\
 &  u_j(x)\to Ax 
 \text{ in } L^1(B_1)\Bigr\} \nonumber
\end{align}
where the set of piecewise constant functions $PC(B_1,\R^N)$ was defined after (\ref{eqdefpsirelpc}).
We remark that periodicity does not play a role here, as $g$ will be used for the local energy
in a representative volume element, where the macroscopic slip field is approximately affine with gradient $A$.
{The minimization in (\ref{eqdefg})
is by itself a complex problem, which in general can only be attacked numerically. The study of efficient tools for this procedure
and of possible explicit solutions in special cases
constitutes an interesting direction of further work.}

%We remark that $\psi_0^\rel(sb,t)=\infty$ whenever $sb\not\in \Z^M$, therefore $g_\infty$ is infinite also on many rank-one matrices.
It is important to notice that
the function $g$ is quasiconvex, continuous, one homogeneous, and has linear growth, in the sense that 
 $\frac1c |A|\le g(A)\le c|A|$ for all matrices $A$.
A detailed proof of those facts will appear elsewhere \cite{ContiGarroniMuellerBumi}.

%\section{Convergence result}
%\label{sec:result}
\def\netw{\mathrm{netw}}
We now pass from the local picture to a macroscopic slip field.
Let $u:\T\to\R^N$ be a slip field, which for now we assume to be continuously differentiable. 
Around every point $x\in \T$, the local energy of the dislocations corresponding to $\nabla u$ is 
determined by the energy $g(\nabla u)$ of the optimal network realizing this gradient, 
as defined in (\ref{eqdefg}). We therefore obtain the expression
\begin{equation}
E_\netw[u]= \int_\T g(\nabla u) dx\,,
\end{equation}
valid for $u\in C^1(\T;\R^N)$.
If the slip field is not smooth, an additional approximation procedure is needed. This leads to the expression
\begin{equation}\label{eqdefeself}
 E_\self[u]=\int_\T g(\nabla u) dx + \int_\T g(\frac{dD^su}{d|D^su|})d|D^su|
\end{equation}
for $u\in BV(\T;\R^N)$. Here $D^su$ is the singular part of the distributional gradient of $u$, which may contain
both jump parts and diffuse parts, and $\nabla u$ is the absolutely continuous part.
Inded, for any $u\in BV(\T;\R^N)$ one writes $Du=\nabla u \calL^2+D^su$, where $\nabla u\in L^1(\T;\R^{N\times 2})$ 
and $D^su$ is a measure on $\T$ with values in $\R^{N\times 2}$, concentrated on a set of Lebesgue measure equal to zero.
Typical { instances} of singular measures are a Dirac measure concentrated on a segment, corresponding to the
derivative of a discontinuous function, and a measure concentrated on a fractal, such as the Cantor set.
{ For example, for a given matrix  $A_*\in \R^{N\times 2}$, the slip
\begin{equation}
 u_*(x)=
 \begin{cases}
  A_*x &\text{ if } 0<x_1<L/2\\
0 & \text{ otherwise}
  \end{cases}
\end{equation}
would have a continuous dislocation density $A_*$ for $0<x_1<L/2$, and a concentrated dislocation density
along the $\{x_1=0\}$ and $\{x_1=L/2\}$ lines. 
The jump in $u_*$ at the point $(0,t)$ equals $u_*(0^+,t)-u_*(0^-,t)=A_*e_2t$; the distributional gradient
is oriented along the normal to the singularity, which is $e_1$. The same holds for $x_1=L/2$.
Inserting in the above expression we obtain that the self-energy associated with $u_*$ reads
\begin{equation}%\label{eqdefeself}
 E_\self[u_*]=\int_{\{0<x_1<L/2\}} g(A_*) dx + 2\int_{0}^L g(tA_*e_2\otimes e_1)dt\,.
\end{equation}}% 
We refer to \cite{AmbrosioFP} for further details.
The fact that $E_\self$ is the relaxation of $E_\netw$ follows, for example, from the general results
on relaxation of $BV$ functionals in  \cite{AmbrosioDalmaso1992} or \cite{FonsecaMueller1993}.
In particular, this implies that $E_\self$ is lower semicontinuous.

\section{The effective limiting energy}
\label{sec:result}
The limiting process discussed above, with subsequent rescalings, can be directly obtained starting from the
phase-field model in a suitable scaling regime which permits the total length of the dislocation lines to diverge as $\eps\to0$.
{ As discussed above, we focus on the regime 
in which the elastic displacement scales as $\ln \frac1\eps$, so that the relevant convergence of 
displacements is  (see discussion around (\ref{eqconvu}))
\begin{equation}\label{eqconvu2}
%\lim_{\eps\to0}
\frac{u_\eps}{\ln \frac1\eps}\to u_0\,.
\end{equation}
Correspondingly, the energy $E_\eps$ is assumed to be proportional to $(\ln \frac1\eps)^2$. We will determine
the limit (in the sense of $\Gamma$-convergence) of the quotients
\begin{equation}
 \frac{1}{(\ln \frac1\eps)^2} E_\eps \stackrel{\Gamma}{\to} E_0
\end{equation}
as in (\ref{eqconveeps}). 
This scaling was chosen so that the self-interaction leading to line-tension effects
balances the macroscopic elastic energy.}

The limiting functional is  given by the sum of the self-energy and the elastic energy,
\begin{equation}\label{eqdefF}
E_0[u]= E_\self[u]+\int_{\T\times\T} K(z) (u(x)-u(x+z))\cdot(u(x)-u(x+z))\,dx\,dz\,,
\end{equation}
with $E_\self$ defined in (\ref{eqdefeself}).
The proof, which is mathematically very technical and builds upon $BV$ relaxation techniques from  \cite{AmbrosioDalmaso1992} and \cite{FonsecaMueller1993} and the lower bound for the case of finite dislocation length in \cite{ContiGarroniMueller2011}, will 
be presented elsewhere  \cite{ContiGarroniMuellerBumi}.
One key idea is that the self-energy $E_\self$ arises from the short-range part of the interaction. Precisely,  one can show
that, for any 
$\rho>0$, 
\begin{equation}
\begin{split}
\frac{1}{(\ln \frac1\eps)^2} 
 \int_\T \left(\int_{B_\rho(0)} K(z) (u_\eps(x)-u_\eps(x+z))\cdot(u_\eps(x)-u_\eps(x+z)) dz \right)dx \\
 \to E_\self[u]\,,
 \end{split}
\end{equation}
for an optimal sequence $u_\eps$ converging to $u$ as in (\ref{eqconvu2}). 
The fact that the limit does not depend on $\rho$ corresponds to the fact that the line-tension
energy $\psi_0$, which then originates $\psi_\rel$, $\psi_\infty$ and $g$, is localized on 
a small neighbourhood of the dislocation, whose size is, as $\eps\to0$, much larger than $\eps$, but smaller than the (fixed)
length $\rho$.

At the same time, the part of the energy
with $|z|>\rho$ is continuous in the limit, and
\begin{align}
\nonumber
\frac{1}{(\ln \frac1\eps)^2} 
 \int_\T \left(\int_{\T\setminus B_\rho(0)} K(z) (u_\eps(x)-u_\eps(x+z))\cdot(u_\eps(x)-u_\eps(x+z)) dz \right)dx \\
 \to
 \int_\T \left(\int_{\T\setminus B_\rho(0)} K(z) (u(x)-u(x+z))\cdot(u(x)-u(x+z)) dz \right)dx \,.
\end{align}
The proof is based on making these two assertions rigorous,  showing that the optimal
sequence can be chosen to be the same for both, and taking the limit $\rho\to0$.

We finally transfer our result back into the 3D setting discussed in the introduction.
Let $U:\T\times \R\to\R^3$ be a deformation, with slip field $\gamma=[U]$. 
The self-energy of the slip field is then obtained using
the change of variables performed in (\ref{eqgammau}) to express the slip field
in a basis of Burgers vectors. Precisely, one obtains
\begin{equation}
 F_\self[\gamma]=\int_\T f(\nabla \gamma) dx + \int_\T f(\frac{dD^s\gamma}{d|D^s\gamma|})d|D^s\gamma|
\end{equation}
for $\gamma\in BV(\T;\R^2\times\{0\})$, where $f$ is defined by the change
of variables 
$f(\sum_{i=1}^N s_i \otimes c_i)=g(\sum_{i=1}^N e_i \otimes c_i)$ for any $c_i\in \R^2$
(here $\{s_i\}$ denotes the basis of $\calB$ chosen in  (\ref{eqgammau})
and $\{e_i\}$ the canonical basis of $\R^N$).
Therefore the functional
\begin{equation}%\label{eqdefhate}
%\frac1{(\ln (1/\eps)^2} F_\eps[U]=
\frac1{(\ln (1/\eps)^2} \left[\frac1\eps\int_\T W_\calB(\gamma_\eps) dx+ \int_{\T\times \R} \frac12 \C e(U_\eps)\cdot e(U_\eps) \, dx\right]
\end{equation}
(with $\gamma_\eps=[U_\eps]$) $\Gamma$-converges, as $\eps\to0$, to the functional
\begin{equation}\label{eqdeffzero}
F_0[U]=\int_\T f(\nabla \gamma) dx + \int_\T f(\frac{dD^s\gamma}{d|D^s\gamma|})d|D^s\gamma|
+
\int_{\T\times \R} \frac12 \C e(U)\cdot e(U) \, dx
\end{equation}
where $\gamma=[U]$. The relevant convergence of the displacement fields is 
$U_\eps/\ln (1/\eps)\to U$, which implies 
$\gamma_\eps/\ln (1/\eps)\to \gamma$ for the corresponding slip fields.
{
We remark that, expressed in the physical dislocation density tensor $\alpha\sim u_0 \ln \frac1\eps$, the
self-energy of the networks scales as $(\ln \frac1\eps)^2F_0[u_0] \sim \ln \frac1\eps \alpha$. }

\section{Conclusions and outlook}
We have discussed a new facet of the phase-field model developed by Michael Ortiz and coworkers in the 90s,
and shown that it permits a rigorous analytical study of dislocation networks in the plane, 
which accounts for the self-energy of the dislocations and the long-range elastic energy.
We discussed how the  effective energy for  dislocation microstructures can be computed by solving
an appropriate cell problem in each representative volume element.
The analysis of the cell problem itself leads to the prediction, under specific circumstances, 
of complex dislocation networks in the plane, with interaction between dislocations in different directions, as illustrated in Figure \ref{fig-micro1}.
Our derivation of the  functional $F_0$ supports the use of  strain-gradient models in plasticity with
 linear growth of the strain gradient term, as was done 
for example in \cite{ContiOrtiz05,AnguigeDondl2014}.
At the same time the current derivation adds more details to these models, as the self-energy is rigorously
derived from the material properties. In particular, the cell structure illustrated in 
Figure \ref{fig-micro1} is a direct consequence of the specific form of the self-energy in cubic crystals, 
and would not appear with the simpler self-energy used in 
\cite{ContiOrtiz05}.
{ Although the  framework we formulated is completely general, our specific analysis is restricted to a cubic 
geometry.  It would be interesting to extend the analysis  to the commonly seen dislocation structures, as for example
to the geometry appropriate for $111$ planes of fcc crystals, which have a triangular structure.
}

\section*{Acknowledgements}
We are very grateful to Michael Ortiz for sharing with us his inspiration and many of his insights
on dislocations, crystal plasticity and much more.
SC and SM acknowledge financial support by the Deutsche Forschungsgemeinschaft through the Sonderforschungsbereich 1060 {\sl ``The mathematics of emergent effects''}, project A5. AG acknowledges the financial support of the Italian National Research Project PRIN 2010 (Calcolo delle Variazioni), 2010A2TFX2$\_$003.

%\bibliographystyle{amsplain-initials}
%\bibliographystyle{alpha-noname}
%\bibliography{cogaor}

\section*{References}


\begin{thebibliography}{49}
\expandafter\ifx\csname natexlab\endcsname\relax\def\natexlab#1{#1}\fi
\providecommand{\url}[1]{\texttt{#1}}
\providecommand{\href}[2]{#2}
\providecommand{\path}[1]{#1}
\providecommand{\DOIprefix}{doi:}
\providecommand{\ArXivprefix}{arXiv:}
\providecommand{\URLprefix}{URL: }
\providecommand{\Pubmedprefix}{pmid:}
\providecommand{\doi}[1]{\href{http://dx.doi.org/#1}{\path{#1}}}
\providecommand{\Pubmed}[1]{\href{pmid:#1}{\path{#1}}}
\providecommand{\bibinfo}[2]{#2}
\ifx\xfnm\relax \def\xfnm[#1]{\unskip,\space#1}\fi
%Type = Article
\bibitem[{Ortiz(1999)}]{Ortiz1999}
\bibinfo{author}{M.~Ortiz},
\newblock \bibinfo{title}{Plastic yielding as a phase transition},
\newblock \bibinfo{journal}{J. Appl. Mech.-Trans. ASME} \bibinfo{volume}{66}
  (\bibinfo{year}{1999}) \bibinfo{pages}{289--298}.
%Type = Article
\bibitem[{Koslowski et~al.(2002)Koslowski, Cuiti{\~n}o, and
  Ortiz}]{KoslowskiCuitinoOrtiz2002}
\bibinfo{author}{M.~Koslowski}, \bibinfo{author}{A.~M. Cuiti{\~n}o},
  \bibinfo{author}{M.~Ortiz},
\newblock \bibinfo{title}{A phase-field theory of dislocation dynamics, strain
  hardening and hysteresis in ductile single crystal},
\newblock \bibinfo{journal}{J. Mech. Phys. Solids} \bibinfo{volume}{50}
  (\bibinfo{year}{2002}) \bibinfo{pages}{2597--2635}.
%Type = Article
\bibitem[{Koslowski and Ortiz(2004)}]{KoslowskiOrtiz2004}
\bibinfo{author}{M.~Koslowski}, \bibinfo{author}{M.~Ortiz},
\newblock \bibinfo{title}{A multi-phase field model of planar dislocation
  networks},
\newblock \bibinfo{journal}{Model. Simul. Mat. Sci. Eng.} \bibinfo{volume}{12}
  (\bibinfo{year}{2004}) \bibinfo{pages}{1087--1097}.
%Type = Book
\bibitem[{Hirth and Lothe(1968)}]{HirthLothe1968}
\bibinfo{author}{J.~P. Hirth}, \bibinfo{author}{J.~Lothe},
  \bibinfo{title}{{Theory of Dislocations}}, \bibinfo{publisher}{McGraw-Hill},
  \bibinfo{address}{New York}, \bibinfo{year}{1968}.
%Type = Book
\bibitem[{Hull and Bacon(2011)}]{HullBacon}
\bibinfo{author}{D.~Hull}, \bibinfo{author}{D.~J. Bacon},
  \bibinfo{title}{Introduction to dislocations}, \bibinfo{edition}{5th} ed.,
  \bibinfo{publisher}{Butterworth-Heinemann}, \bibinfo{address}{Oxford, UK},
  \bibinfo{year}{2011}.
%Type = Article
\bibitem[{Xu and Ortiz(1993)}]{XuOrtiz1993}
\bibinfo{author}{G.~Xu}, \bibinfo{author}{M.~Ortiz},
\newblock \bibinfo{title}{A variational boundary integral method for the
  analysis of {3-D} cracks of arbitrary geometry modelled as continuous
  distributions of dislocation loops},
\newblock \bibinfo{journal}{International journal for numerical methods in
  engineering} \bibinfo{volume}{36} (\bibinfo{year}{1993})
  \bibinfo{pages}{3675--3701}.
%Type = Article
\bibitem[{Xu and Argon(2000)}]{XuArgon2000}
\bibinfo{author}{G.~Xu}, \bibinfo{author}{A.~S. Argon},
\newblock \bibinfo{title}{Homogeneous nucleation of dislocation loops under
  stress in perfect crystals},
\newblock \bibinfo{journal}{Philosophical magazine letters}
  \bibinfo{volume}{80} (\bibinfo{year}{2000}) \bibinfo{pages}{605--611}.
%Type = Article
\bibitem[{Xiang et~al.(2008)Xiang, Wei, Ming, and Weinan}]{XiangWeiMingE2008}
\bibinfo{author}{Y.~Xiang}, \bibinfo{author}{H.~Wei},
  \bibinfo{author}{P.~Ming}, \bibinfo{author}{E.~Weinan},
\newblock \bibinfo{title}{A generalized {Peierls}--{Nabarro} model for curved
  dislocations and core structures of dislocation loops in {Al} and {Cu}},
\newblock \bibinfo{journal}{Acta materialia} \bibinfo{volume}{56}
  (\bibinfo{year}{2008}) \bibinfo{pages}{1447--1460}.
%Type = Article
\bibitem[{Garroni and M{\"u}ller(2005)}]{GarroniMueller2005}
\bibinfo{author}{A.~Garroni}, \bibinfo{author}{S.~M{\"u}ller},
\newblock \bibinfo{title}{{$\Gamma$}-limit of a phase-field model of
  dislocations},
\newblock \bibinfo{journal}{SIAM J. Math. Anal.} \bibinfo{volume}{36}
  (\bibinfo{year}{2005}) \bibinfo{pages}{1943--1964}.
%Type = Article
\bibitem[{Garroni and M{\"u}ller(2006)}]{GarroniMueller2006}
\bibinfo{author}{A.~Garroni}, \bibinfo{author}{S.~M{\"u}ller},
\newblock \bibinfo{title}{A variational model for dislocations in the line
  tension limit},
\newblock \bibinfo{journal}{Arch. Ration. Mech. Anal.} \bibinfo{volume}{181}
  (\bibinfo{year}{2006}) \bibinfo{pages}{535--578}.
%Type = Article
\bibitem[{Cacace and Garroni(2009)}]{CacaceGarroni2009}
\bibinfo{author}{S.~Cacace}, \bibinfo{author}{A.~Garroni},
\newblock \bibinfo{title}{A multi-phase transition model for the dislocations
  with interfacial microstructure},
\newblock \bibinfo{journal}{Interfaces Free Bound.} \bibinfo{volume}{11}
  (\bibinfo{year}{2009}) \bibinfo{pages}{291--316}.
%Type = Article
\bibitem[{Conti et~al.(2011)Conti, Garroni, and
  M\"uller}]{ContiGarroniMueller2011}
\bibinfo{author}{S.~Conti}, \bibinfo{author}{A.~Garroni},
  \bibinfo{author}{S.~M\"uller},
\newblock \bibinfo{title}{Singular kernels, multiscale decomposition of
  microstructure, and dislocation models},
\newblock \bibinfo{journal}{Arch. Rat. Mech. Anal.} \bibinfo{volume}{199}
  (\bibinfo{year}{2011}) \bibinfo{pages}{779--819}.
%Type = Article
\bibitem[{Conti and Gladbach(2015)}]{ContiGladbach2015}
\bibinfo{author}{S.~Conti}, \bibinfo{author}{P.~Gladbach},
\newblock \bibinfo{title}{A line-tension model of dislocation networks on
  several slip planes},
\newblock \bibinfo{journal}{Mechanics of Materials} \bibinfo{volume}{90}
  (\bibinfo{year}{2015}) \bibinfo{pages}{140--147}.
%Type = Article
\bibitem[{Conti et~al.(2015)Conti, Garroni, and
  Massaccesi}]{ContiGarroniMassaccesi2015}
\bibinfo{author}{S.~Conti}, \bibinfo{author}{A.~Garroni},
  \bibinfo{author}{A.~Massaccesi},
\newblock \bibinfo{title}{Modeling of dislocations and relaxation of
  functionals on 1-currents with discrete multiplicity},
\newblock \bibinfo{journal}{Calc. Var. PDE} \bibinfo{volume}{54}
  (\bibinfo{year}{2015}) \bibinfo{pages}{1847--1874}.
%Type = Article
\bibitem[{Ariza and Ortiz(2005)}]{ArizaOrtiz06}
\bibinfo{author}{M.~P. Ariza}, \bibinfo{author}{M.~Ortiz},
\newblock \bibinfo{title}{Discrete crystal plasticity},
\newblock \bibinfo{journal}{Arch. Ration. Mech. Anal.} \bibinfo{volume}{178}
  (\bibinfo{year}{2005}) \bibinfo{pages}{149--226}.
%Type = Article
\bibitem[{Ramasubramaniam et~al.(2007)Ramasubramaniam, Ariza, and
  Ortiz}]{RamaArizaOrtiz2007}
\bibinfo{author}{A.~Ramasubramaniam}, \bibinfo{author}{M.~P. Ariza},
  \bibinfo{author}{M.~Ortiz},
\newblock \bibinfo{title}{A discrete mechanics approach to dislocation dynamics
  in bcc crystals},
\newblock \bibinfo{journal}{J. Mech. Phys. Solids} \bibinfo{volume}{55}
  (\bibinfo{year}{2007}) \bibinfo{pages}{615--647}.
%Type = Article
\bibitem[{Conti et~al.(2015{\natexlab{a}})Conti, Garroni, and
  Ortiz}]{ContiGarroniOrtiz2015}
\bibinfo{author}{S.~Conti}, \bibinfo{author}{A.~Garroni},
  \bibinfo{author}{M.~Ortiz},
\newblock \bibinfo{title}{The line-tension approximation as the dilute limit of
  linear-elastic dislocations},
\newblock \bibinfo{journal}{Arch. Ration. Mech. Anal.} \bibinfo{volume}{218}
  (\bibinfo{year}{2015}{\natexlab{a}}) \bibinfo{pages}{699--755}.
%Type = Article
\bibitem[{Conti et~al.(2015{\natexlab{b}})Conti, Garroni, and
  Ortiz}]{ContiGarroniOrtizDiscrete}
\bibinfo{author}{S.~Conti}, \bibinfo{author}{A.~Garroni},
  \bibinfo{author}{M.~Ortiz},
\newblock \bibinfo{title}{From atomic interactions to the line-tension
  approximation for dilute dislocations},
\newblock \bibinfo{journal}{in preparation}
  (\bibinfo{year}{2015}{\natexlab{b}}).
%Type = Article
\bibitem[{Barnett and Swanger(1971)}]{BarnettSwanger1971}
\bibinfo{author}{D.~Barnett}, \bibinfo{author}{L.~Swanger},
\newblock \bibinfo{title}{The elastic energy of a straight dislocation in an
  infinite anisotropic elastic medium},
\newblock \bibinfo{journal}{Physica Status Solidi} \bibinfo{volume}{48}
  (\bibinfo{year}{1971}) \bibinfo{pages}{419--428}.
%Type = Incollection
\bibitem[{Rice(1985)}]{Rice:1985b}
\bibinfo{author}{J.~R. Rice},
\newblock \bibinfo{title}{Conserved integrals and energetic forces},
\newblock in: \bibinfo{editor}{K.~J. Miller} (Ed.),
  \bibinfo{booktitle}{Fundamentals of Deformation and Fracture},
  \bibinfo{publisher}{Cambridge University Press}, \bibinfo{year}{1985}.
%Type = Book
\bibitem[{Gottstein(2014)}]{Gottstein2014materialwissenschaft}
\bibinfo{author}{G.~Gottstein}, \bibinfo{title}{Materialwissenschaft und
  Werkstofftechnik: Physikalische Grundlagen},
  \bibinfo{publisher}{Springer-Verlag}, \bibinfo{year}{2014}.
%Type = Article
\bibitem[{Fleck and Hutchinson(1993)}]{FleckHutchinson1993}
\bibinfo{author}{N.~A. Fleck}, \bibinfo{author}{J.~W. Hutchinson},
\newblock \bibinfo{title}{A phenomenological theory for strain gradient effects
  in plasticity},
\newblock \bibinfo{journal}{J. Mech. Phys. Solids} \bibinfo{volume}{41}
  (\bibinfo{year}{1993}) \bibinfo{pages}{1825--1857}.
%Type = Article
\bibitem[{Nix and Gao(1998)}]{NixGao1998}
\bibinfo{author}{W.~D. Nix}, \bibinfo{author}{H.~J. Gao},
\newblock \bibinfo{title}{Indentation size effects in crystalline materials: A
  law for strain gradient plasticity},
\newblock \bibinfo{journal}{J. Mech. Phys. Solids} \bibinfo{volume}{46}
  (\bibinfo{year}{1998}) \bibinfo{pages}{411--425}.
%Type = Article
\bibitem[{Fleck and Hutchinson(2001)}]{FleckHutchinson2001}
\bibinfo{author}{N.~A. Fleck}, \bibinfo{author}{J.~W. Hutchinson},
\newblock \bibinfo{title}{A reformulation of strain gradient plasticity},
\newblock \bibinfo{journal}{J. Mech. Phys. Solids} \bibinfo{volume}{49}
  (\bibinfo{year}{2001}) \bibinfo{pages}{2245--2271}.
%Type = Article
\bibitem[{Bassani(2001)}]{Bassani2001}
\bibinfo{author}{J.~L. Bassani},
\newblock \bibinfo{title}{Incompatibility and a simple gradient theory of
  plasticity},
\newblock \bibinfo{journal}{J. Mech. Phys. Solids} \bibinfo{volume}{49}
  (\bibinfo{year}{2001}) \bibinfo{pages}{1983--1996}.
%Type = Article
\bibitem[{Conti and Ortiz(2005)}]{ContiOrtiz05}
\bibinfo{author}{S.~Conti}, \bibinfo{author}{M.~Ortiz},
\newblock \bibinfo{title}{Dislocation microstructures and the effective
  behavior of single crystals},
\newblock \bibinfo{journal}{Arch. Rat. Mech. Anal.} \bibinfo{volume}{176}
  (\bibinfo{year}{2005}) \bibinfo{pages}{103--147}.
%Type = Article
\bibitem[{Kuroda and Tvergaard(2008)}]{KurodaTvergaard2008a}
\bibinfo{author}{M.~Kuroda}, \bibinfo{author}{V.~Tvergaard},
\newblock \bibinfo{title}{On the formulations of higher-order strain gradient
  crystal plasticity models},
\newblock \bibinfo{journal}{J. Mech. Phys. Solids} \bibinfo{volume}{56}
  (\bibinfo{year}{2008}) \bibinfo{pages}{1591--1608}.
%Type = Article
\bibitem[{Fokoua et~al.(2014)Fokoua, Conti, and Ortiz}]{FokouaContiOrtiz2014}
\bibinfo{author}{L.~Fokoua}, \bibinfo{author}{S.~Conti},
  \bibinfo{author}{M.~Ortiz},
\newblock \bibinfo{title}{Optimal scaling laws for ductile fracture derived
  from strain-gradient microplasticity},
\newblock \bibinfo{journal}{J. Mech. Phys. Solids} \bibinfo{volume}{62}
  (\bibinfo{year}{2014}) \bibinfo{pages}{295--311}.
%Type = Article
\bibitem[{Aubry and Ortiz(2003)}]{AubryOrtiz03}
\bibinfo{author}{S.~Aubry}, \bibinfo{author}{M.~Ortiz},
\newblock \bibinfo{title}{The mechanics of deformation-induced
  subgrain-dislocation structures in metallic crystals at large strains},
\newblock \bibinfo{journal}{Proc. R. Soc. Lond. A} \bibinfo{volume}{459}
  (\bibinfo{year}{2003}) \bibinfo{pages}{3131--3158}.
%Type = Article
\bibitem[{Ortiz and Repetto(1999)}]{OrtizRepetto1999}
\bibinfo{author}{M.~Ortiz}, \bibinfo{author}{E.~Repetto},
\newblock \bibinfo{title}{Nonconvex energy minimization and dislocation
  structures in ductile single crystals},
\newblock \bibinfo{journal}{J. Mech. Phys. Solids} \bibinfo{volume}{47}
  (\bibinfo{year}{1999}) \bibinfo{pages}{397--462}.
%Type = Article
\bibitem[{Conti and Theil(2005)}]{ContiTheil2005}
\bibinfo{author}{S.~Conti}, \bibinfo{author}{F.~Theil},
\newblock \bibinfo{title}{Single-slip elastoplastic microstructures},
\newblock \bibinfo{journal}{Arch. Rat. Mech. Anal.} \bibinfo{volume}{178}
  (\bibinfo{year}{2005}) \bibinfo{pages}{125--148}.
%Type = Article
\bibitem[{Conti et~al.(2013)Conti, Dolzmann, and
  Kreisbeck}]{ContiDolzmannKreisbeck2013b}
\bibinfo{author}{S.~Conti}, \bibinfo{author}{G.~Dolzmann},
  \bibinfo{author}{C.~Kreisbeck},
\newblock \bibinfo{title}{Relaxation of a model in finite plasticity with two
  slip systems},
\newblock \bibinfo{journal}{Math. Models Methods Appl. Sci.}
  \bibinfo{volume}{23} (\bibinfo{year}{2013}) \bibinfo{pages}{2111--2128}.
%Type = Article
\bibitem[{Anguige and Dondl(2014)}]{AnguigeDondl2014}
\bibinfo{author}{K.~Anguige}, \bibinfo{author}{P.~Dondl},
\newblock \bibinfo{title}{Relaxation of the single-slip condition in
  strain-gradient plasticity},
\newblock \bibinfo{journal}{R. Soc. Lond. Proc. Ser. A Math. Phys. Eng. Sci.}
  \bibinfo{volume}{470} (\bibinfo{year}{2014}).
%Type = Article
\bibitem[{Hunter and Koslowski(2008)}]{HunterKoslowski2008}
\bibinfo{author}{A.~Hunter}, \bibinfo{author}{M.~Koslowski},
\newblock \bibinfo{title}{Direct calculations of material parameters for
  gradient plasticity},
\newblock \bibinfo{journal}{Journal of the Mechanics and Physics of Solids}
  \bibinfo{volume}{56} (\bibinfo{year}{2008}) \bibinfo{pages}{3181--3190}.
%Type = Article
\bibitem[{Xiang(2009)}]{Xiang2009}
\bibinfo{author}{Y.~Xiang},
\newblock \bibinfo{title}{Continuum approximation of the peach--koehler force
  on dislocations in a slip plane},
\newblock \bibinfo{journal}{Journal of the Mechanics and Physics of Solids}
  \bibinfo{volume}{57} (\bibinfo{year}{2009}) \bibinfo{pages}{728--743}.
%Type = Article
\bibitem[{Zhu and Xiang(2014)}]{ZhuXiang2014}
\bibinfo{author}{X.~Zhu}, \bibinfo{author}{Y.~Xiang},
\newblock \bibinfo{title}{Continuum framework for dislocation structure, energy
  and dynamics of dislocation arrays and low angle grain boundaries},
\newblock \bibinfo{journal}{Journal of the Mechanics and Physics of Solids}
  \bibinfo{volume}{69} (\bibinfo{year}{2014}) \bibinfo{pages}{175--194}.
%Type = Article
\bibitem[{Garroni et~al.(2010)Garroni, Leoni, and
  Ponsiglione}]{GarroniLeoniPonsiglione2010}
\bibinfo{author}{A.~Garroni}, \bibinfo{author}{G.~Leoni},
  \bibinfo{author}{M.~Ponsiglione},
\newblock \bibinfo{title}{Gradient theory for plasticity via homogenization of
  discrete dislocations},
\newblock \bibinfo{journal}{J. Eur. Math. Soc. (JEMS)} \bibinfo{volume}{12}
  (\bibinfo{year}{2010}) \bibinfo{pages}{1231--1266}.
%Type = Article
\bibitem[{M\"uller et~al.(2014)M\"uller, Scardia, and
  Zeppieri}]{MuellerScardiaZeppieri}
\bibinfo{author}{S.~M\"uller}, \bibinfo{author}{L.~Scardia},
  \bibinfo{author}{C.~I. Zeppieri},
\newblock \bibinfo{title}{Geometric rigidity for incompatible fields and an
  application to strain-gradient plasticity},
\newblock \bibinfo{journal}{Indiana Univ. Math. J.} \bibinfo{volume}{63}
  (\bibinfo{year}{2014}) \bibinfo{pages}{1365--1396}.
%Type = Incollection
\bibitem[{M{\"u}ller et~al.(2015)M{\"u}ller, Scardia, and
  Zeppieri}]{MuellerScardiaZeppieri2015}
\bibinfo{author}{S.~M{\"u}ller}, \bibinfo{author}{L.~Scardia},
  \bibinfo{author}{C.~I. Zeppieri},
\newblock \bibinfo{title}{Gradient theory for geometrically nonlinear
  plasticity via the homogenization of dislocations},
\newblock in: \bibinfo{editor}{S.~Conti}, \bibinfo{editor}{K.~Hackl} (Eds.),
  \bibinfo{booktitle}{Analysis and Computation of Microstructure in Finite
  Plasticity}, \bibinfo{publisher}{Springer}, \bibinfo{year}{2015}, pp.
  \bibinfo{pages}{175--204}.
%Type = Article
\bibitem[{Ambrosio and Braides(1990{\natexlab{a}})}]{AmbrosioBraides1990a}
\bibinfo{author}{L.~Ambrosio}, \bibinfo{author}{A.~Braides},
\newblock \bibinfo{title}{Functionals defined on partitions in sets of finite
  perimeter. {I}.\ {I}ntegral representation and {$\Gamma$}-convergence},
\newblock \bibinfo{journal}{J. Math. Pures Appl. (9)} \bibinfo{volume}{69}
  (\bibinfo{year}{1990}{\natexlab{a}}) \bibinfo{pages}{285--305}.
%Type = Article
\bibitem[{Ambrosio and Braides(1990{\natexlab{b}})}]{AmbrosioBraides1990b}
\bibinfo{author}{L.~Ambrosio}, \bibinfo{author}{A.~Braides},
\newblock \bibinfo{title}{Functionals defined on partitions in sets of finite
  perimeter. {II}.\ {S}emicontinuity, relaxation and homogenization},
\newblock \bibinfo{journal}{J. Math. Pures Appl. (9)} \bibinfo{volume}{69}
  (\bibinfo{year}{1990}{\natexlab{b}}) \bibinfo{pages}{307--333}.
%Type = Book
\bibitem[{Ambrosio et~al.(2000)Ambrosio, Fusco, and Pallara}]{AmbrosioFP}
\bibinfo{author}{L.~Ambrosio}, \bibinfo{author}{N.~Fusco},
  \bibinfo{author}{D.~Pallara}, \bibinfo{title}{Functions of Bounded Variation
  and Free Discontinuity Problems}, Mathematical Monographs,
  \bibinfo{publisher}{Oxford University Press}, \bibinfo{year}{2000}.
%Type = Book
\bibitem[{{Dal Maso}(1993)}]{Dalmaso1993}
\bibinfo{author}{G.~{Dal Maso}}, \bibinfo{title}{An introduction to
  {$\Gamma$}-convergence}, Progress in Nonlinear Differential Equations and
  their Applications, 8, \bibinfo{publisher}{Birkh\"auser Boston Inc.},
  \bibinfo{address}{Boston, MA}, \bibinfo{year}{1993}.
%Type = Book
\bibitem[{Braides(2002)}]{Braides02}
\bibinfo{author}{A.~Braides}, \bibinfo{title}{Gamma-convergence for Beginners},
  \bibinfo{publisher}{Oxford University Press}, \bibinfo{year}{2002}.
%Type = Article
\bibitem[{Bulatov and Kaxiras(1997)}]{BlatovKaxiras1997}
\bibinfo{author}{V.~V. Bulatov}, \bibinfo{author}{E.~Kaxiras},
\newblock \bibinfo{title}{Semidiscrete variational peierls framework for
  dislocation core properties},
\newblock \bibinfo{journal}{Physical review letters} \bibinfo{volume}{78}
  (\bibinfo{year}{1997}) \bibinfo{pages}{4221}.
%Type = Article
\bibitem[{Dai et~al.(2013)Dai, Xiang, and Srolovitz}]{DaiXiangSrolovitz2013}
\bibinfo{author}{S.~Dai}, \bibinfo{author}{Y.~Xiang}, \bibinfo{author}{D.~J.
  Srolovitz},
\newblock \bibinfo{title}{Structure and energy of (111) low-angle twist
  boundaries in {Al}, {Cu} and {Ni}},
\newblock \bibinfo{journal}{Acta materialia} \bibinfo{volume}{61}
  (\bibinfo{year}{2013}) \bibinfo{pages}{1327--1337}.
%Type = Article
\bibitem[{Conti et~al.(2015)Conti, Garroni, and
  M\"uller}]{ContiGarroniMuellerBumi}
\bibinfo{author}{S.~Conti}, \bibinfo{author}{A.~Garroni},
  \bibinfo{author}{S.~M\"uller},
\newblock \bibinfo{title}{Derivation of the dislocation self-energy from a
  planar phase-field model},
\newblock \bibinfo{journal}{in preparation}  (\bibinfo{year}{2015}).
%Type = Article
\bibitem[{Ambrosio and Dal~Maso(1992)}]{AmbrosioDalmaso1992}
\bibinfo{author}{L.~Ambrosio}, \bibinfo{author}{G.~Dal~Maso},
\newblock \bibinfo{title}{On the relaxation in {${\rm BV}(\Omega;{\bf R}^m)$}
  of quasi-convex integrals},
\newblock \bibinfo{journal}{J. Funct. Anal.} \bibinfo{volume}{109}
  (\bibinfo{year}{1992}) \bibinfo{pages}{76--97}.
%Type = Article
\bibitem[{Fonseca and M{\"u}ller(1993)}]{FonsecaMueller1993}
\bibinfo{author}{I.~Fonseca}, \bibinfo{author}{S.~M{\"u}ller},
\newblock \bibinfo{title}{Relaxation of quasiconvex functionals in {${\rm
  BV}(\Omega,{\bf R}^p)$} for integrands {$f(x,u,\nabla u)$}},
\newblock \bibinfo{journal}{Arch. Rational Mech. Anal.} \bibinfo{volume}{123}
  (\bibinfo{year}{1993}) \bibinfo{pages}{1--49}.

\end{thebibliography}
\end{document}